\documentclass{amsart}

\newtheorem{theorem}{Theorem}
\newtheorem{lemma}[theorem]{Lemma}

\newtheorem{corollary}[theorem]{Corollary}

\theoremstyle{definition}

\usepackage{amscd,amssymb}

\begin{document}

\title[Trace cocharacters of two $3\times 3$ matrices]
{Multiplicities in the mixed trace cocharacter sequence of two
$3\times 3$ matrices}

\author[Drensky, Genov, and Valenti]
{Vesselin Drensky, Georgi K. Genov, and Angela Valenti}
\address{Institute of Mathematics and Informatics,
Bulgarian Academy of Sciences, 1113 Sofia, Bulgaria}
\email{drensky@math.bas.bg, gguenov@hotmail.com}
\address{Dipartimento di Matematica ed Applicazioni, Universit\`a di Palermo,
Via Archirafi 34, 90123 Palermo, Italy}
\email{avalenti@math.unipa.it}
\thanks{The research of the first two authors was partially supported by Grant
MM-1106/2001 of the Bulgarian Foundation for Scientific Research.}
\thanks{The research of the third author was partially supported by MIUR, Italy.}

\begin{abstract}
We find explicitly the multiplicities in the (mixed) trace
cocharacter sequence of two $3\times 3$ matrices over a field of
characteristic 0 and show that asymptotically they behave as
polynomials of seventh degree. As a consequence we obtain also the
multiplicities of certain irreducible characters in the
cocharacter sequence of the polynomial identities of $3\times 3$
matrices.
\end{abstract}
\maketitle

\section*{Introduction}

All considerations in this paper are over an arbitrary field $F$
of characteristic 0. Let $n\geq 2$ be a fixed integer. Consider
the $d$ generic $n\times n$ matrices $X_1,\ldots,X_d$, $d\geq 2$.
The main results of our paper concern the case $n=3$ and $d=2$.
There are several important algebras related to $X_1,\ldots,X_d$.
Among them are the algebra $R_d$ generated by $X_1,\ldots,X_d$,
the pure (or commutative) trace algebra $C_d$ generated by the
traces of all products $\text{tr}(X_{i_1}\cdots X_{i_k})$, and the
mixed (or noncommutative) trace algebra $T_d$ generated by $R_d$
and $C_d$ regarding the elements of $C_d$ as scalar matrices. We
denote by $C,R$, and $T$ the corresponding algebras related to a
countable set $\{X_1,X_2,\ldots\}$ of generic matrices.

The algebra $R$ is one of the most important objects in the theory
of algebras with polynomial identities. It is isomorphic to the
factor algebra $F\langle x_1,x_2,\ldots\rangle/I(M_n(F))$ of the
free associative algebra $F\langle x_1,x_2,\ldots\rangle$ modulo
the ideal $I(M_n(F))$ of the polynomial identities of the $n\times
n$ matrix algebra $M_n(F)$. The algebra $C_d$ has a natural
interpretation in classical invariant theory, as the algebra of
invariants of the general linear group $GL_n(F)$ acting by
simultaneous conjugation on $d$ matrices of size $n$. The algebra
$T_d$ is known as the algebra of matrix concominants and also
consists of the invariant functions under a suitable action of
$GL_n(F)$. See e.g. the books \cite{J}, \cite{F3}, or \cite{DF} as
a background on $R_d,C_d$, and $T_d$ and their application to
invariant theory, structure theory of PI-algebras, and theory of
finite dimensional division algebras.

The algebras $R_d,C_d,T_d$ as well as $R,C,T$ are graded by
multidegree. The symmetric group $S_k$ acts naturally on the
multilinear elements of degree $k$ of $R,T,C$. The corresponding
$S_k$-characters $\chi_k(R)=\chi_k(M_n(F)),\chi_k(C),\chi_k(T)$
are called, respectively, the cocharacter of the polynomial
identities, the pure trace cocharacter, and the mixed trace
cocharacter of $n\times n$ matrices. They decompose as
\[
\chi_k(M_n(F))=\sum_{\lambda\vdash
k}m_{\lambda}(M_n(F))\chi_{\lambda},
\]
\[
\chi_k(C)=\sum_{\lambda\vdash k}m_{\lambda}(C)\chi_{\lambda},
\]
\[
\chi_k(T)=\sum_{\lambda\vdash k}m_{\lambda}(T)\chi_{\lambda},
\]
where $\lambda$ is a partition of $k$ and $\chi_{\lambda}$ is the
related irreducible $S_k$-character. If some of the multiplicities
$m_{\lambda}(M_n(F)), m_{\lambda}(C), m_{\lambda}(T)$ is nonzero,
then $\lambda=(\lambda_1,\ldots,\lambda_{n^2})$ is a partition in
not more than $n^2$ parts.

The algebras $R_d,C_d,T_d$ are also $GL_d(F)$-modules with
$GL_d(F)$-action induced by the canonical $GL_d(F)$-action on the
vector space with basis $\{X_1,\ldots,X_d\}$ and decompose as
\[
R_d=\sum_{k\geq 0}\sum_{\lambda\vdash
k}m_{\lambda}(M_n(F))W_{\lambda}, \quad C_d=\sum_{k\geq
0}\sum_{\lambda\vdash k}m_{\lambda}(C)W_{\lambda},\quad
T_d=\sum_{k\geq 0}\sum_{\lambda\vdash k}m_{\lambda}(T)W_{\lambda},
\]
where $W_{\lambda}$ is the irreducible $GL_d(F)$-module related to
$\lambda=(\lambda_1,\ldots,\lambda_d)$, and the multiplicities of
$W_{\lambda}$ are the same as the multiplicities of
$\chi_{\lambda}$ in the corresponding $S_k$-cocharacter.

The Hilbert (or Poincar\'e) series of $R_d$ is
\[
H(R_d)=H(R_d,t_1,\ldots,t_d)=\sum\dim R_d^{(k)}t_1^{k_1}\cdots
t_d^{k_d},
\]
where $R_d^{(k)}$ is the homogeneous component of multidegree
$k=(k_1,\ldots,k_d)$. Similarly one defines the Hilbert series of
$C_d$ and $T_d$. These series are symmetric functions and
decompose as infinite linear combinations of Schur functions
$S_{\lambda}(t_1,\ldots,t_d)$. Since the Hilbert series play the
role of characters of $GL_d(F)$ and the Schur functions are the
characters of the corresponding irreducible $GL_d(F)$-modules, the
multiplicities $m_{\lambda}(M_n(F))$, $m_{\lambda}(C)$, and
$m_{\lambda}(T)$ can be obtained from the Hilbert series of
$R_d,C_d$, and $T_d$ for $d=n^2$. For small $d$, the Hilbert
series give the multiplicities for
$\lambda=(\lambda_1,\ldots,\lambda_d)$ only. See the book by
Macdonald \cite{M} as a standard reading on theory of symmetric
functions and the book by one of the authors \cite{D} for the
applications of representation theory of $S_k$ and $GL_d(F)$ to
PI-algebras.

The algebras $C$ and $T$ are described in terms of invariant
theory and are easier to study than the algebra $R$. Since the
three algebras are very close to each other, the standard way to
investigate the polynomial identities of $M_n(F)$ is via $C$ and
$T$. The multiplicities of the cocharacters of $M_n(F),C,T$ are
explicitly found for $n=2$ only, by Formanek \cite{F1} for
$M_2(F),C,T$, and, with different methods, by Drensky \cite{D1}
for $M_2(F)$, see also Procesi \cite{P}. In particular,
\[
m_{\lambda}(T)=(\lambda_1-\lambda_2+1)(\lambda_2-\lambda_3+1)(\lambda_3-\lambda_4+1),
\]
if $\lambda=(\lambda_1,\lambda_2,\lambda_3,\lambda_4)$ and
$m_{\lambda}(T)=0$ if $\lambda_5\not=0$.

It follows from classical invariant theory that the Hilbert series
of $C_d$ and $T_d$ can be expressed as multiple integrals but for
$n\geq 3$ their direct evaluation is quite difficult and was given
by Teranishi \cite{T1, T2} for $C_2$ and $n=3$ and $n=4$ only. Van
den Bergh \cite{V} found a graph theoretical approach for the
calculation of $H(C_d)$ and $H(T_d)$. Berele and Stembridge
\cite{BS} calculated the Hilbert series of $C_d$ and $T_d$ for
$n=3$, $d\leq 3$ and of $T_2$ for $n=4$, correcting also some
typographical errors in the expression of $H(C_2)$ for $n=4$ in
\cite{T2}. Using the Hilbert series of $C_2$, $n=3$, Berele
\cite{B} found an asymptotic expression of
$m_{(\lambda_1,\lambda_2)}(C)$. The explicit form of the
generating function of $m_{(\lambda_1,\lambda_2)}(C)$ was found by
Drensky and Genov \cite{DG1} correcting also a technical error in
\cite{B}. The proof of \cite{DG1} is quite technical. Later
Drensky and Genov \cite{DG2} suggested a method to find the
coefficients of the Schur functions in the expansion of a class of
rational symmetric functions in two variables. This significantly
simplified their proof of the formula for the generating function
of $m_{(\lambda_1,\lambda_2)}(C)$. A recent general result of
Carbonara, Carini and Remmel \cite{CCR} states that, for any $n$,
there exists a positive integer $u_n$, such that for a fixed
partition $\bar\lambda=(\lambda_2,\ldots,\lambda_{n^2})$,  there
exist $u_n$ polynomials $P^{\bar\lambda}_{n,k}(\lambda_1)$,
$k=0,1,\ldots,u_n-1$, of degree $n-1$ and with the same leading
term, with the property that
$m_{\lambda}(C)=P^{\bar\lambda}_{n,k}(\lambda_1)$, for $\lambda_1$
sufficiently large, where
$\lambda=(\lambda_1,\bar\lambda)=(\lambda_1,\ldots,\lambda_{n^2})$,
$k\equiv\lambda_1 \mod u_n$. Since the cocharacter sequences of
$C$ and $T$ are related by $\chi_k(T)=\chi_{k+1}(C)\downarrow
S_k$, and the multiplicities  $m_{\lambda}(T)$ can be expressed in
terms of the multiplicities of $C$ by the branching theorem, this
easily implies a similar statement for $m_{\lambda}(T)$. For $n=3$
the expression of $u_n$ given in \cite{CCR} is $u_3=6$.

The purpose of the present paper is to calculate, for $n=3$, the
generating function of $m_{(\lambda_1,\lambda_2)}(T)$ and to find
explicit expressions for $m_{(\lambda_1,\lambda_2)}(T)$. As in the
case of $m_{(\lambda_1,\lambda_2)}(C)$, see \cite{B}, the
multiplicity $m_{(\lambda_1,\lambda_2)}(T)$ behaves asymptotically
as a polynomial of degree 7 in $\lambda_1$ and $\lambda_2$. Our
approach is to apply the method of Drensky and Genov \cite{DG2} to
the explicit form of the Hilbert series of $T_2$ found by Berele
and Stembridge \cite{BS}. As a consequence of some results of
Formanek \cite{F1, F2} we obtain the exact values of the
multiplicities $m_{\lambda}(M_3(F))$ for
$\lambda=(\lambda_1,\ldots,\lambda_9)$ when
$\lambda_3=\cdots=\lambda_9\geq 2$. The calculations have been
performed independently by the first and the third authors in
Palermo and by the second author in Sofia, as an additional
warranty of their correctness. In both cases we have used Maple.

\section{Preliminaries}
Let $F[[x,y]]$ be the algebra of formal power series in two
variables and let
\[
\text{Sym}[[x,y]]=F[[x,y]]^{S_2}
\]
be the subalgebra of symmetric functions. The set of Schur
functions
\[
\left\{S_{\lambda}(x,y)\mid \lambda=(\lambda_1,\lambda_2)
\in{\mathbb Z}\times {\mathbb Z}, \ \lambda_1\geq\lambda_2\geq 0
\right\}
\]
forms a basis of $\text{Sym}[[x,y]]$ as a topological vector space
with its usual formal power series topology.  Hence for any
symmetric function $ f(x,y)$ in two variables,
\[
f(x,y)=\sum_{i,j\geq 0}a(i,j)x^iy^j
=\sum_{\lambda}m(\lambda)S_{\lambda}(x,y),
\]
where $a(i,j)=a(j,i)\in F$ for all $i,j\geq 0$ and $m(\lambda)\in
F$ is the multiplicity of $S_{\lambda}(x,y)$ in the decomposition
of $f(x,y)$. The Schur function $S_{\lambda}(x,y)$ has the
following simple expression
\[
S_{\lambda}(x,y)=(xy)^{\lambda_2}(x^p+x^{p-1}y+\cdots+xy^{p-1}+y^p)=
\frac{(xy)^{\lambda_2}\left(x^{p+1}-y^{p+1}\right)}{x-y}, \eqno(1)
\]
where we have denoted $p=\lambda_1-\lambda_2$. This easily gives
that the coefficients $a(i,j)$ and the multiplicities $m(\lambda)$
of $f(x,y)$  are related by
\[
m(\lambda_1,\lambda_2)=a(\lambda_1,\lambda_2)-a(\lambda_1+1,\lambda_2-1),
\]
where $a(\lambda_1+1,\lambda_2-1)=0$ if $\lambda_2=0$.

Drensky and Genov \cite{DG1} introduced the multiplicity series of
$f(x,y)$, namely
\[
M(f)(t,u)=\sum_{\lambda_1\geq\lambda_2\geq 0}
m(\lambda_1,\lambda_2)t^{\lambda_1}u^{\lambda_2} \in F[[t,u]].
\]
Introducing a new variable $v=tu$, the series $M(f)(t,u)$ accepts
the more convenient form
\[
M'(f)(t,v)=M(f)(t,u)=\sum_{\lambda_1\geq\lambda_2\geq 0}
m(\lambda_1,\lambda_2)t^{\lambda_1-\lambda_2}v^{\lambda_2} \in
F[[t,v]]
\]
and the mapping $M':\text{Sym}[[x,y]]\to F[[t,v]]$ is a continuous
linear bijection.

The relation between symmetric functions and their multiplicity
series is given, see (11) in \cite{DG1}, by
\[
f(x,y)=\frac{xM'(f)(x,xy)-yM'(f)(y,xy)}{x-y}. \eqno(2)
\]
The symmetric functions depending on the product $xy$ only behave
like constants under the mapping $M'$. If $a(z)\in F[[z]]$ is a
formal power series in one variable $z$, then
\[
M'(a(xy)f(x,y))=a(v)M'(f(x,y)), \eqno(3)
\]
for any symmetric function $f(x,y)$, see (10) in \cite{DG1}.

We consider two linear operators $Y$ and $Y_a$ of $F[[t,v]]$. If
$h(t,v)\in F[[t,v]]$, there exists a unique $f(x,y)\in
\text{Sym}[[x,y]]$ such that $h(t,v)=M'(f)$. As in \cite{DG1}, we
define
\[
Y(h)(t,v)=M'\left(\frac{f(x,y)}{(1-x)(1-y)}\right).
\]
Similarly, if $a(z)\in F[[z]]$, then
\[
Y_a(h)(t,v)=M'\left(\frac{f(x,y)}{(1-a(xy)x)(1-a(xy)y)}\right).
\]
Clearly, $Y=Y_1$. The next lemma gives relations between
$h(t,v),Y(h)$, and $Y_a(h)$.

\begin{lemma}\label{lemma 1}
If $h(t,v)\in F[[t,v]]$ and $a(z)\in F[[z]]$, then
\[
Y(h)(t,v)=\frac{th(t,v)-vh(v,v)}{(1-t)(t-v)}, \eqno(4)
\]
\[
Y_a(h)(t,v)=\frac{th(t,v)-a(v)vh(a(v)v,v)}{(1-a(v)t)(t-a(v)v)}.
\eqno(5)
\]
\end{lemma}
\begin{proof} The first equation is established in Proposition 2
of \cite{DG1} and the second in Lemma 1 of \cite{DG2}. Since
\cite{DG2} is an announcement only, we give the proof for
completeness of the exposition. Since the mapping $M'$ is linear
and continuous, it is sufficient to verify (5) for the image
$h(t,v)=M'(S_{\lambda}(x,y))$ of the Schur functions only. By (1)
and (3) we may assume that $\lambda_2=0$, i.e. $\lambda=(p,0)$,
and we obtain that $h(t,v)=t^p$. Direct verification shows that
the right hand side of (5) has the form
\[
h_a(t,v)=\frac{t^{p+1}-(a(v)v)^{p+1}}{(1-a(v)t)(t-a(v)v)}.
\]
Further calculations show that
\[
xh_a(x,xy)=\frac{x^{p+1}-(xy)^{p+1}a(xy)^{p+1}}{(1-a(xy)x)(1-a(xy)y)},
\]
\[
\frac{xh_a(x,xy)-yh_a(y,xy)}{x-y}=\frac{x^{p+1}-y^{p+1}}{(x-y)(1-a(xy)x)(1-a(xy)y)}
\]
\[
=\frac{S_{(p,0)}(x,y)}{(1-a(xy)x)(1-a(xy)y)},
\]
i.e. $S_{(p,0)}(x,y)$ and $h_a(t,v)$ satisfy the relation (2).
Since $M'$ is a bijection, this completes the proof.
\end{proof}

Finally, we need the following equality which is a translation in
our language of a partial case of a result of Thrall \cite{Th},
see also (8) in \cite{DG1}:
\[
M'\left(\frac{1}{(1-x^2)(1-xy)(1-y^2)}\right)
=\frac{1}{(1-t^2)(1-v^2)}.\eqno(6)
\]

\section{Main Results}
Till the end of the paper we consider $3\times 3$ matrices only,
i.e. we fix $n=3$. Most of our considerations are for $d=2$ and we
replace the variables $t_1,t_2$ with $x,y$. The Hilbert series of
$T_2$ found by Berele and Stembridge \cite{BS} is
\[
H(T_2,x,y)=\frac{1}{(1-x)^2(1-y)^2
(1-x^2)(1-y^2)(1-xy)^2(1-x^2y)(1-xy^2)}.\eqno(7)
\]

The first of our main results gives the multiplicity series of
this Hilbert series.

\begin{theorem}
The multiplicity series of the Hilbert series $H(T_2,x,y)$ of the
mixed trace algebra $T_2$ of two generic $3\times 3$ matrices is
\[
M'(H(T_2,x,y))(t,v)=\frac{h_3(v)t^3+h_2(v)t^2+h_1(v)t+h_0(v)}
{(1-v)^7(1+v)^4(1+v^2)(1-t)^3(1+t)(1-vt)},\eqno(8)
\]
where the polynomials $h_i(v)\in F[v]$, $i=0,1,2,3$, are
\[
h_3(v)=v^2(v^4-v^3+3v^2-v+1),\quad h_2(v)= v(2v^4-4v^3+v^2-v-1),
\]
\[
h_1(v)=v(-v^4-v^3+v^2-4v+2), \quad h_0(v)=v^4-v^3+3v^2-v+1.
\]
It has also the expression
\[
M'(H(T_2,x,y))(t,v)=\frac{a_3(v)}{(1-t)^3}+\frac{a_2(v)}{(1-t)^2}
+\frac{a_1(v)}{(1-t)}+\frac{b(v)}{1+t}+\frac{c(v)}{1-vt},\eqno(9)
\]
where
\[
a_3(v)=\frac{1}{2(1-v)^6(1+v)^2},\quad
a_2(v)=\frac{(3v^2-2v+1)}{2^2(1-v)^7(1+v)^3},
\]
\[
a_1(v)=\frac{(v^4-6v^3+14v^2-6v+1)} {2^3(1-v)^8(1+v)^4},
\]
\[
b(v)=\frac{1}{2^3(1-v)^2(1+v)^4(1+v^2)}, \quad
c(v)=\frac{-v^4}{(1-v)^8(1+v)^4(1+v^2)}.
\]
\end{theorem}

\begin{proof}
Direct calculations, which we have performed using Maple, show
that the expressions from the right hand sides of (8) and (9) are
equal. So, the proof will be completed, if we show that the
replacement of any of these expressions in (2) gives the Hilbert
series in (7). This has been done also using Maple.

We think that it is interesting to know how we have calculated (8)
and (9). We rewrite $H(T_2,x,y)$ in the form
\[
\frac{1}{1-xy}\left(\frac{1}{(1-x)(1-y)}\right)^2
\frac{1}{(1-(xy)x)(1-(xy)y)} \left(\frac{1}
{(1-x^2)(1-xy)(1-y^2)}\right).
\]
It follows from (3), (4), and (5) that
\[
M'(H(T_2,x,y))(t,v)
=\frac{1}{1-v}Y^2Y_{z}M'\left(\frac{1}{(1-x^2)(1-xy)(1-y^2)}\right).
\]
Taking into account (6), we obtain
\[
M'(H(T_2,x,y))(t,v)
=\frac{1}{1-v}Y^2Y_{z}\left(\frac{1}{(1-t^2)(1-v^2)}\right).
\]
Now we calculate consecutively
\[
w_0(t,v)=\frac{1}{(1-t^2)(1-v^2)},
\]
\[
w_1=Y_z(w_0) =\frac{1+v^2t}{(1-v^2)^2(1+v^2)(1-t^2)(1-vt)},
\]
\[
w_2=Y(w_1)=\frac{(-v^2(v^2-v+1)t^2-v(v^2-1)t+(v^2-v+1))}{(1-v)^4(1+v)^3(1+v^2)(1-t)^2(1+t)(1-vt)},
\]
\[
w_3=Y(w_2)=\frac{h_3(v)t^3+h_2(v)t^2+h_1(v)t+h_0(v)}
{(1-v)^6(1+v)^4(1+v^2)(1-t)^3(1+t)(1-vt)},
\]
\[
M'(H(T_2,x,y))(t,v)=w_4=\frac{w_3}{1-v},
\]
as in (8). The expression (9) is obtained considering (8) as a
rational function in $t$ with coefficients from $F(v)$ and
presenting it as a sum of elementary fractions.
\end{proof}

\begin{lemma}\label{Lemma 3} The function $M'(H(T_2,x,y))(t,v)$
can be presented in the form
\[
\sum_{p,q \geq 0}\left(a^+_{pq} + (-1)^qa^-_{pq} + (-1)^pb^+_q +
(-1)^{p+q}b^-_q \right)t^pv^q -\frac{1}{64}\sum_{p,r \geq 0}
(-1)^{p+r}t^pv^{2r+1}
\]
\[
+\sum_{p,s\geq 0}\left[c^+_s + (-1)^sc^-_s\right](tv)^pv^s
+\frac{1}{64}\sum_{p,w \geq 0}(-1)^w(tv)^pv^{2w},
\]
where
\[
a^+_{pq}
=\frac{(84p^2+14pq+q^2)q^5}{2^57!}+\frac{(40p^2+12pq+q^2)q^4}{2^55!}
+\frac{(90p^2+49pq+5q^2)q^3}{2^73^2}
\]
\[
+\frac{(104p^2+108pq+15q^2)q^2}{2^73} +\frac{(19596p^2+43666pq
+9599q^2)q}{2^66!}
\]
\[
+\frac{1800p^2+11676pq+4993q^2}{2^65!}+\frac{9492p+11437q}{2^9\cdot
3 \cdot 7}+\frac{43}{64},
\]
\[
a^-_{pq}= \frac{(12p^2+18pq+7q^2)q}{2^{10}3} +
\frac{8p^2+28pq+17q^2}{2^9} +\frac{180p+229q}{2^93}+\frac{13}{64},
\]
\[
b^+_q = \frac{q+4}{2^8},\quad b^-_q =
\frac{2q^3+24q^2+85q+84}{2^83},
\]
\[
c^+_s= -\frac{s^7}{2^5 7!}-\frac{s^6}{2^6 5!} -\frac{19s^5}{2^56!}
-\frac{s^4}{2^8 3!} + \frac{391s^3}{2^6 6!}+ \frac{79s^2}{2^{10}5}
-\frac{1453s}{2^7 5!7}-\frac{17}{2^{10}},
\]
\[
c^-_s = -\frac{s^3+9s^2+17s-3}{2^{10}3}.
\]
\end{lemma}
\begin{proof} We decompose the rational functions
$a_3(v),a_2(v),a_1(v),b(v),c(v)\in F(v)$ as linear combinations of
elementary fractions of the form
\[
\frac{1}{(1-v)^k},\quad \frac{1}{(1+v)^k},\quad \frac{1}{1+v^2},
\quad \frac{v}{1+v^2}.
\]
The results are
\[
a_3(v)=\frac{1}{8(1-v)^6}+\frac{1}{8(1-v)^5}+\frac{3}{32(1-v)^4}
+\frac{1}{16(1-v)^3}
\]
\[
+\frac{5}{128(1-v)^2}+\frac{3}{128(1-v)}+\frac{1}{128(1+v)^2}+\frac{3}{128(1+v)},
\]
\[
a_2(v)=\frac{1}{16(1-v)^7}-\frac{1}{32(1-v)^6}+\frac{1}{32(1-v)^4}
+\frac{11}{256(1-v)^3}
\]
\[
+\frac{21}{512(1-v)^2}+\frac{17}{512(1-v)}+\frac{3}{256(1+v)^3}
+\frac{13}{512(1+v)^2}+\frac{17}{512(1+v)},
\]
\[
a_1(v)=\frac{1}{32(1-v)^8}-\frac{1}{32(1-v)^6}-\frac{1}{32(1-v)^5}-\frac{5}{512(1-v)^4}
\]
\[
+\frac{3}{256(1-v)^3}+\frac{25}{1024(1-v)^2}+\frac{29}{1024(1-v)}
\]
\[
+\frac{7}{512(1+v)^4}+\frac{7}{256(1+v)^3}
+\frac{33}{1024(1+v)^2}+\frac{29}{1024(1+v)},
\]
\[
b(v)=\frac{1}{256(1-v)^2}+\frac{3}{256(1-v)}+\frac{1}{64(1+v)^4}+\frac{1}{32(1+v)^3}
\]
\[
+\frac{9}{256(1+v)^2}+\frac{7}{256(1+v)}-\frac{v}{64(1+v^2)},
\]
\[
c(v)=-\frac{1}{32(1-v)^8}+\frac{1}{32(1-v)^7}+\frac{1}{32(1-v)^6}-\frac{11}{512(1-v)^4}
\]
\[
-\frac{11}{512(1-v)^3}-\frac{9}{1024(1-v)^2}+\frac{1}{256(1-v)}
\]
\[
-\frac{1}{512(1+v)^4}-\frac{1}{512(1+v)^3}
+\frac{1}{1024(1+v)^2}+\frac{1}{256(1+v)}+\frac{1}{64(1+v^2)}.
\]
Applying the formula
\[
\frac{1}{(1-z)^{k+1}}=\sum_{m\geq 0}\binom{k+m}{k}z^m\eqno(10)
\]
for $z=\pm v, -v^2, \pm t, vt$, and expressing the binomial
coefficients in (10) as polynomials of degree $k$ in $m$, we
obtain the presentation of $M'(H(T_2,x,y))(t,v)$ in the statement
of the lemma.
\end{proof}

The following theorem gives explicit formulas for the
multiplicities $m_{(\lambda_1,\lambda_2)}(T)$.

\begin{theorem} The multiplicities $m_{(\lambda_1,\lambda_2)}(T)$
of the mixed trace cocharacter of $3\times 3$ matrices are given
by the following formulas, where we have presented
$(\lambda_1,\lambda_2)$ in the form $(p+q,q)$:
\[
m_{(\lambda_1,\lambda_2)}(T)=a^+_{pq} + (-1)^qa^-_{pq} +
(-1)^pb^+_q + (-1)^{p+q}b^-_q -\frac{(-1)^{p+r}\varepsilon_1}{64}
\]
\[
+\delta\left(c^+_{q-p}+(-1)^{q-p}c^-_{q-p}+\frac{(-1)^{p+w}\varepsilon_2}{64}\right).
\]
Here $\varepsilon_1=1$, if $q=2r+1$ and $\varepsilon_1=0$, if
$q=2r$; $\varepsilon_2=1$, if $q-p=2w$ and $\varepsilon_2=0$, if
$q-p=2w+1$; and $\delta=0$, if $\lambda_1>2\lambda_2$ and
$\delta=1$, if $\lambda_1\leq 2\lambda_2$. The expressions of
$a^{\pm}_{pq},b^{\pm}_q,c^{\pm}_{s}$ are given in the above Lemma
3.
\end{theorem}

\begin{proof}
Clearly, the multiplicity
$m_{(\lambda_1,\lambda_2)}(T)=m_{(p+q,q)}(T)$ is equal to the
coefficient of $t^pv^q$ in the expansion of the formal power
series $M'(H(T_2,x,y))(t,v)$. This explains the contribution of
$a^{\pm}_{pq},b^{\pm}_q$ and  $\varepsilon_1/64$ to
$m_{(p+q,q)}(T)$. The coefficient $c^{\pm}_s$ contributes to
$(tv)^pv^s$, which is equal to $t^pv^q$ for $s=q-p$. Hence we need
$q\geq p$ which is equivalent to $\lambda_1\leq 2\lambda_2$. In
this case $\delta=1$. Otherwise $\delta=0$. The coefficient
$\varepsilon_2$ appears by the same reasons as $\varepsilon_1$.
\end{proof}

The expression of $m_{(\lambda_1,\lambda_2)}(T)$  becomes much
simpler if we are interested in their asymptotic behavior only.
The following corollary follows immediately from Theorem 4,
presenting $\lambda$ in the form $(p+q,q)$. It is in the spirit of
the description of $m_{(\lambda_1,\lambda_2)}(C)$ given by Berele
\cite{B}.

\begin{corollary}\label{corollary 5}
For any partition $\lambda=(\lambda_1,\lambda_2)$, the
multiplicities $m_{\lambda}=m_{\lambda}(T)$ of the mixed trace
cocharacter of $3\times 3$ matrices satisfy the condition
\[
m_{\lambda}=\frac{\lambda_2^7}{7!2^5}+\frac{(\lambda_1-\lambda_2)\lambda_2^6}{6!2^4}
+\frac{(\lambda_1-\lambda_2)^2\lambda_2^5}{5!2^4} +{\mathcal
O}((\lambda_1+\lambda_2)^6),
\]
if $\lambda_1>2\lambda_2\geq 0$ and
\[
m_{\lambda}=\frac{\lambda_2^7}{7!2^5}+\frac{(\lambda_1-\lambda_2)\lambda_2^6}{6!2^4}
+\frac{(\lambda_1-\lambda_2)^2\lambda_2^5}{5!2^4}
-\frac{(2\lambda_2-\lambda_1)^7}{7!2^5} +{\mathcal
O}((\lambda_1+\lambda_2)^6),
\]
if $2\lambda_2\geq \lambda_1\geq \lambda_2\geq 0$.

If $\lambda_2$ is fixed, then, for $\lambda_1>2\lambda_2$,
\[
m_{\lambda}=\lambda_1^2\left(\left(\frac{\lambda_2^5}{5!2^4}
+\frac{\lambda_2^4}{2^5\cdot
3}+\frac{5\lambda_2^3}{2^6}+\frac{13\lambda_2^2}{4!2}+\frac{1633\lambda_2}{5!2^5}+\frac{15}{2^6}\right)
\right.
\]
\[
\left.
+(-1)^{\lambda_2}\left(\frac{\lambda_2}{2^8}+\frac{1}{2^6}\right)\right)+
{\mathcal O}(\lambda_1).
\]
\end{corollary}

Comparing with the asymptotic expression of
$m_{(\lambda_1,\lambda_2)}(C)$ in the form of Theorem 15 in
\cite{DG1}, we see that
\[
m_{(\lambda_1,\lambda_2)}(C)\approx
\frac{1}{9}m_{(\lambda_1,\lambda_2)}(T).
\]
The formulas for $m_{(\lambda_1,\lambda_2)}(C)$ in \cite{DG1}
agree with the result of \cite{CCR} that, for a fixed
$\lambda_2,\ldots,\lambda_9$, and for $\lambda_1$ sufficiently
large, the multiplicity $m_{\lambda}(C)$ behaves, depending on
$\lambda_1 \mod 6$, like a polynomial of second degree in
$\lambda_1$, with a leading term which does not depend on
$\lambda_1 \mod 6$. The second part of Corollary \ref{corollary 5}
is in the same spirit. The careful study of the form of
$M'(H(T_2,x,y))$ from Theorem 2 shows that it is sufficient to
consider $\lambda_1$ mod 2 and not mod 6.

Finally, we give an application to the ``ordinary'' cocharacters
of the polynomial identities of $3\times 3$ matrices.

\begin{corollary}
If $\mu=(\mu_1,\ldots,\mu_9)$ is such that $\mu_3=\cdots=\mu_9\geq
2$, then the multiplicity $m_{\mu}(M_3(F))$ is equal to the
multiplicity $m_{(\lambda_1,\lambda_2)}(T)$ found in this paper,
where $\lambda_1=\mu_1-\mu_3$, $\lambda_2=\mu_2-\mu_3$.
\end{corollary}

\begin{proof}
By a result of Formanek \cite{F1}, $m_{\mu}(T)=m_{\lambda}(T),$
for any partitions $\lambda=(\lambda_1,\ldots,\lambda_{n^2})$ and
$\mu=(\mu_1,\ldots,\mu_{n^2})$ related by the equalities
$\mu_k=\lambda_k+r$ for all $k=1,\ldots,n^2$ and all $n\geq 2$.
Also, for $\lambda_{n^2}$ sufficiently large,
$m_{\lambda}(M_n(F))=m_{\lambda}(T)$. Another result of Formanek,
in \cite{F2}, gives that the latter equality holds already for
$\lambda_{n^2}\geq 2$. This immediately completes the proof.
\end{proof}

\section*{Acknowledgements}
A large part of this project was carried out when the first author
visited the University of Palermo with the financial support of
the INdAM, Italy. He is very grateful for the hospitality and the
creative atmosphere during his stay in Palermo.


\begin{thebibliography}{99}
\bibitem{B}
A. Berele, {\it Approximate multiplicities in the trace
cocharacter sequence of two three-by-three matrices}, Commun.
Algebra {\bf 25} (1997), 1975-1983.
\bibitem{BS}
A. Berele, J.R. Stembridge, {\it Denominators for the Poincar\'e
series of invariants of small matrices}, Israel J. Math. {\bf 114}
(1999), 157-175.
\bibitem{CCR}
J.O. Carbonara, L. Carini, J.B. Remmel, {\it Trace cocharacters
and the Kronecker products of Schur functions}, J. Algebra {\bf
260} (2003), 631-656.
\bibitem{D1}
V. Drensky, {\it Codimensions of T-ideals and Hilbert series of
relatively free algebras}, J. Algebra {\bf 91} (1984), 1-17.
\bibitem{D}
V. Drensky, {\it Free Algebras and PI-Algebras}, Springer-Verlag,
Singapore, 1999.
\bibitem{DF}
V. Drensky, E. Formanek, {\it Polynomial Identity Rings}, Advanced
Courses in Mathematics, CRM Barcelona, Birkh\"auser Verlag, Basel,
2004.
\bibitem{DG1}
V. Drensky, G.K. Genov, {\it Multiplicities of Schur functions in
invariants of two $3\times 3$ matrices}, J. Algebra {\bf 264}
(2003), 496-519.
\bibitem{DG2}
V. Drensky, G.K. Genov, {\it Multiplicities of Schur functions
with applications to invariant theory and PI-algebras}, C.R. Acad.
Bulg. Sci. {\bf 57} (2004), No. 3, 5-10.
\bibitem{F1}
E. Formanek, {\it Invariants and the ring of generic matrices}, J.
Algebra {\bf 89} (1984), 178-223.
\bibitem{F2}
E. Formanek, {\it A conjecture of Regev about the Capelli
polynomial}, J. Algebra {\bf 109} (1987), 93-114.
\bibitem{F3}
E. Formanek, {\it The Polynomial Identities and Invariants of $n
\times n$ Matrices}, CBMS Regional Conf. Series in Math. {\bf 78},
Published for the Confer. Board of the Math. Sci. Washington DC,
AMS, Providence RI, 1991.
\bibitem{J}
N. Jacobson, {\it PI-Algebras. An Introduction}, Lecture Notes in
Mathematics, {\bf 441}, Springer-Verlag, Berlin-New York, 1975.
\bibitem{M}
I.G. Macdonald, {\it Symmetric Functions and Hall Polynomials},
Oxford Univ. Press (Clarendon), Oxford, 1979, Second Edition,
1995.
\bibitem{P}
C. Procesi, {\it Computing with $2\times 2$ matrices}, J. Algebra
{\bf 87} (1984), 342-359.
\bibitem{T1}
Y. Teranishi, {\it The ring of invariants
of matrices}, Nagoya Math. J. {\bf 104} (1986), 149-161.
\bibitem{T2}
Y. Teranishi, {\it Linear Diophantine equations and invariant
theory of matrices}, Commutative algebra and combinatorics (Kyoto,
1985), 259-275, Adv. Stud. Pure Math., {\bf 11}, North-Holland,
Amsterdam, 1987. \bibitem{Th} R.M. Thrall, {\it On symmetrized
Kronecker powers and the structure of the free Lie ring}, Trans.
Amer. Math. Soc. {\bf 64} (1942), 371-388.
\bibitem{V}
M. Van den Bergh, {\it Explicit rational forms for the Poincar\"e
series of the trace rings of generic matrices}, Isr. J. Math. {\bf
73} (1991), 17-31.
\end{thebibliography}
\end{document}